\documentclass[titlepage,11pt]{article}
\usepackage{latexsym}
\usepackage{amsfonts}
\usepackage{amsmath}
\newcommand\blackslug{\hbox{\hskip 1pt \vrule width 4pt height 8pt depth 1.5pt
        \hskip 1pt}}
\newcommand\bbox{\hfill \quad \blackslug \bigbreak}

\def\ll{,\ldots,}
\newcommand\vare{\varepsilon}
\DeclareMathOperator{\polylog}{polylog}

\title{Pure pairs. VII. Homogeneous submatrices in $0/1$-matrices with a forbidden submatrix}
\author{Alex Scott
\\
Mathematical Institute, University of Oxford, Oxford OX2 6GG, UK
\\
\\
Paul Seymour\thanks{Supported by AFOSR grant
A9550-19-1-0187 and NSF grant  DMS-1800053.}\\
Princeton University, Princeton, NJ 08544
\\
\\
Sophie Spirkl\thanks{This material is based upon work supported by the National Science
Foundation under Award No. DMS-1802201.}\\
University of Waterloo, Ontario, Canada N2L3G1
}

\date{}

\newtheorem{thm}{}[section]

\newcommand{\Proof}{\noindent{\bf Proof.}\ \ }

\begin{document}
\maketitle
\begin{abstract}
For integer $n>0$, let $f(n)$ be the number of rows of the largest
all-$0$ or all-$1$ square submatrix of $M$, minimized over all $n\times n$ $0/1$-matrices $M$.
Thus $f(n)= O(\log n)$. But let us
fix a matrix $H$, and define $f_H(n)$ to be the same, minimized over over all $n\times n$ $0/1$-matrices $M$ such that neither $M$ 
nor its complement (that is, change all $0$'s to $1$'s
and vice versa) 
contains $H$ as a submatrix. It is known that $f_H(n)\ge \vare n^c$, where $c, \vare>0$ are constants depending on $H$.

When can we take $c=1$? If so, then one of $H$ and its complement must be an acyclic matrix (that is, 
the corresponding bipartite graph is a forest).
Kor\'andi, Pach, and Tomon~\cite{pach} conjectured the converse, that $f_H(n)$ is linear in $n$ for 
every acyclic matrix $H$; and
they proved it for certain matrices $H$ with only two rows.  

Their conjecture remains open, but we show $f_H(n)=n^{1-o(1)}$ for every acyclic matrix~$H$; and indeed
there is a $0/1$-submatrix that is either $\Omega(n)\times n^{1-o(1)}$ or $n^{1-o(1)}\times \Omega(n)$. 
\end{abstract}

\section{Introduction}

A $0/1$-matrix can be regarded as a bipartite graph, with a distinguished bipartition $(V_1,V_2)$ say, in which there are
linear orders imposed on $V_1$ and on $V_2$. Submatrix containment corresponds, in graph theory terms, to induced subgraph
containment, respecting the two bipartitions and preserving the linear orders.
In two earlier papers~\cite{trees,pure4}
(one with with Maria Chudnovsky), we proved some results about excluding induced subgraphs, in a general graph and in
a bipartite graph respectively. Now we impose orders on the vertex sets, and only consider induced subgraph containment that 
respects the orders; and we ask how far our earlier theorems remain true under this much weaker hypothesis. 

In this paper, all graphs are finite and with no loops or parallel edges. 
Two disjoint sets are {\em complete} to each other if every vertex of the first is adjacent to every vertex of the second, and
{\em anticomplete} if there are no edges between them.
A pair $(Z_1,Z_2)$ of subsets of $V(G)$ is {\em pure} if 
$Z_1$ is either complete or anticomplete to $Z_2$. Let us state the earlier theorems that we want to 
extend to ordered graphs. First, we proved the following, with Chudnovsky~\cite{trees}:
\begin{thm}\label{trees}
For every forest $T$, there exists $\vare>0$ such that if $G$ is a graph with $n\ge 2$ vertices, and no induced subgraph
is isomorphic to $T$ or its complement, then there is a pure pair $(Z_1,Z_2)$ of subsets of $V(G)$ with $|Z_1|, |Z_2|\ge \vare n$.
\end{thm}
This theorem characterizes forests: if $T$ is a graph that is not a forest
or the complement of one, then there is no $\vare>0$ as in \ref{trees}.

Second, we proved a similar theorem about bipartite graphs, but for this we need some more definitions. A {\em bigraph} is a graph
together with a bipartition $(V_1(G), V_2(G))$ of $G$. A bigraph $G$ {\em contains} a bigraph $H$ if there is an isomorphism
from $H$ to an induced subgraph of $G$ that maps $V_i(H)$ into $V_i(G)$ for $i = 1,2$. The {\em bicomplement} of
a bigraph $H$ is the bigraph obtained by reversing the adjacency of $v_1,v_2$ for all $v_i\in V_i(G)\;(i=1,2)$.
We proved the following in~\cite{pure4}:
\begin{thm}\label{biptrees}
For every forest bigraph $T$, there exists $\vare>0$ such that if $G$ is a bigraph that
does not contain $T$ or its bicomplement, then there is a pure pair $(Z_1,Z_2)$ with $Z_i\subseteq V_i(G)$ and 
$|Z_i|\ge \vare |V_i(G)|$ for $i = 1,2$.
\end{thm}
Again, this characterizes forests, in that if $H$ is a bigraph that is not a forest or the bicomplement of a forest, there is no
$\vare>0$ as in \ref{biptrees}.

What if we impose an order on the vertex set, and ask for the induced subgraph containment to respect the order?
Let us say an {\em ordered graph} is a graph with a linear order on its vertex set. Every induced subgraph inherits an order on
its vertex set in the natural way: let us say an ordered graph $G$ {\em contains}
an ordered graph $H$ if $H$ is isomorphic to an induced subgraph $H'$ of $G$, where the isomorphism carries the order on
$V(H)$
to the inherited order on $V(H')$. One could ask for an analogue of \ref{trees} for ordered graphs, but it is false.
Fox~\cite{fox} showed:
\begin{thm}\label{fox}
Let $H$ be the three-vertex path with vertices $h_1,h_2,h_3$ in order, and make
$H$ an ordered graph using the same order.
For all sufficiently large $n$, there is an ordered graph $G$ with $n$ vertices, that does not contain $H$,
and such that there do not exist two disjoint subsets of $V(G)$, both of size at least $n/\log(n)$, and complete or
anticomplete.
\end{thm}
To deduce that \ref{trees} does not extend to ordered graphs, let $T$ be an ordered tree such that both $T$ and its bicomplement
contain $H$, and use the construction from \ref{fox}.
But something like \ref{trees} {\em is} true: we proved in~\cite{pure6} that:
\begin{thm}\label{pure6}
For every ordered forest $T$ and all $c>0$, there exists $\vare>0$ such that if $G$ is an ordered graph with $|G|\ge 2$ 
that does not contain $T$ or its
complement, there is a pure pair $(Z_1,Z_2)$ in $G$ with $|Z_1|,|Z_2|\ge \vare|G|^{1-c}$.
\end{thm}

Perhaps the situation is better for ordered bipartite graphs: certainly we are not so well-supplied with counterexamples, 
and there are some positive results about ordered bipartite graphs, proved recently by Kor\'andi, Pach, 
and Tomon~\cite{pach}.
Let us say an {\em ordered bigraph} is a bigraph with linear orders on $V_1(G)$ and on $V_2(G)$. This is just a $0/1$ matrix in disguise,
but graph theory language is convenient for us. (Note that we are not giving a linear
order of $V(G)$: that is much too strong and trivially does not work.) An ordered bigraph $G$ {\em contains} an ordered bigraph $H$
if there is an induced subgraph $H'$ of $G$ and an isomorphism from $H$ to $H'$ mapping $V_i(H)$ to $V_i(H')$
and mapping the order on $V_i(H)$ to the inherited order on $V_i(H')$, for $i = 1,2$. (In matrix language, this is just 
submatrix containment.) Kor\'andi, Pach, and Tomon~\cite{pach} showed:
\begin{thm}\label{2x2}
Let $H$ be an ordered bigraph with $|V_1(H)|\le 2$, such that either
\begin{itemize}
\item  $|V_2(H)|\le 2$ and both $H$ and its bicomplement are forests, or
\item  every vertex in $V_2(H)$ has degree exactly one.
\end{itemize}
Then there exists $\vare>0$ with the following property. Let $G$ be an ordered bigraph that does not contain
$H$, with $|V_1(G)|,|V_2(G)|\ge n$; then there is a pure pair $(Z_1,Z_2)$ with $Z_i\subseteq V_i(G)$ and $|Z_i|\ge \vare n$ for $i = 1,2$.
\end{thm}
In both cases of \ref{2x2}, the bigraph $H$ is a forest and so is its bicomplement. Kor\'andi, Pach, and Tomon asked which other
ordered bigraphs $H$ satisfy the conclusion of \ref{2x2}. They observed that every such bigraph must be a forest and the 
bicomplement of a forest, and conjectured that this was sufficient as well as necessary, that is:
\begin{thm}\label{pachconj}
{\bf Conjecture:} Let $H$ be an ordered bigraph such that both $H$ and its bicomplement are forests.
Then there exists $\vare>0$ with the following property. Let $G$ be an ordered bigraph that does not contain
$H$, with $|V_1(G)|,|V_2(G)|\ge n$; then there is a pure pair $(Z_1,Z_2)$ with $Z_i\subseteq V_i(G)$ and $|Z_i|\ge \vare n$ for $i = 1,2$. 
\end{thm}

We have not been able to decide this conjecture, and indeed
have not even been able to prove it for the forest $H$ consisting of a five-vertex path and an isolated vertex with 
$|V_1(H)|=|V_2(H)|$ (under any ordering of $V_1(H)$ and $V_2(H)$).
But we will prove in \ref{symthm} that, for a much more general class of   
ordered bigraphs, it is possible to find pairs of almost linear size. 

Kor\'andi, Pach, and Tomon also proposed an even stronger conjecture (to see that it implies \ref{pachconj}, let $H$ be 
as in \ref{pachconj}, let $H'$ be an ordered forest that contains both $H$ and its bicomplement, and apply \ref{falseconj} for $H'$):
\begin{thm}\label{falseconj}
{\bf Conjecture:} For every ordered forest bigraph $H$, there exists $\vare>0$ with the following property.
Let $G$ be an ordered bigraph that does not contain
$H$ or its bicomplement, with $|V_1(G)|,|V_2(G)|\ge n$: then 
there is a pure pair $(Z_1,Z_2)$ with $Z_i\subseteq V_i(G)$ and $|Z_i|\ge \vare n$ for $i = 1,2$. 
\end{thm}
This seems a natural extension of \ref{pachconj}, in analogy with \ref{trees} and \ref{biptrees}; but, for what it is worth,
our guess is that \ref{falseconj} is false. Perhaps 
$n/\polylog(n)$ might be true?

There is another result of Kor\'andi, Pach, and Tomon, in the same paper~\cite{pach}:
\begin{thm}\label{stars}
Let $H$ be an ordered forest bigraph such that $|V_1(H)|=2$ and $|V_2(H)|=k$.
For every $\tau >0$, there exists $\delta>0$ with the following property.
Let $G$ be an ordered bigraph that does not contain
$H$, with $|V_1(G)|=|V_2(G)|=n$, and such that its bicomplement has at least $\tau n^2$ edges.
Then there are subsets $Z_i\subseteq V_i(G)$ for $i = 1,2$ 
with $|Z_1|\ge \delta n2^{-(1+o(1))(\log \log (\delta n))^k}$ and $|Z_2|\ge \delta n$,
such that $Z_1,Z_2$ are anticomplete.
\end{thm}

So, $|Z_1|$ is not quite linear, but there is more of significance.
There is nothing here about forbidding $G$ to contain the bicomplement of a forest, since the bicomplement of $H$ 
need not be a forest; and the ``$Z_1$ complete to $Z_2$'' outcome is gone. In compensation they have the assumption that the bicomplement of $G$ is not too sparse.

Our objective in this paper is essentially to generalize \ref{stars} to all ordered forest bigraphs $H$. We will 
give two results. Both prove the existence of 
anticomplete sets $Z_1,Z_2$ of cardinalities at least $n^{1-o(1)}$, but neither implies the other. One result (the second)
gives a linear lower bound for one of the sets and a sublinear bound for the other; and the other result (the first) gives a sublinear (but better) bound for both sets.

Every forest is an induced subgraph of a tree, so
we will assume $H$ is a tree, for convenience. The {\em radius} of a tree $T$ is the minimum $r$ such that for some vertex $v$, every vertex
of $T$ can be joined to $v$ by a path with at most $r$ edges.
In the first half of the paper, we will show:
\begin{thm}\label{mainthm}
Let $T$ be an ordered tree bigraph, of radius $r$, and with $t\ge 2$ vertices.
Let $G$ be a bigraph not containing $T$, with $|V_1(G)|, |V_2(G)|\ge n$, such that every vertex of $G$ has degree at most $n/(4t^2)$.
Choose $K$ such that $t^{K^r}=n$.
Then
there are subsets $Z_i\subseteq V_i(G)$ for $i = 1,2$, with $|Z_1|,|Z_2|\ge nt^{-5K^{r-1}}$, such that $Z_1,Z_2$ are anticomplete.
\end{thm}
In the second half of the paper we will show:
\begin{thm}\label{linmainthm}
Let $T$ be an ordered tree bigraph. For all  $c>0$ there exists $\vare>0$ with the following property.
Let $G$ be a bigraph not containing $T$, such that every vertex in $V_1(G)$ has degree less than $\vare |V_2(G)|$, and 
every vertex in $V_2(G)$ has degree less than $\vare |V_1(G)|$.
Then
there are subsets $Z_i\subseteq V_i(G)$ for $i = 1,2$, either with $|Z_1|\ge \vare |V_1(G)|$ and $|Z_2|\ge \vare |V_2(G)|^{1-c}$, 
or with $|Z_1|\ge \vare |V_1(G)|^{1-c}$ and $|Z_2|\ge \vare |V_2(G)|$, such that $Z_1,Z_2$ are anticomplete.
\end{thm}

The first result, \ref{mainthm}, implies:
\begin{thm}\label{mainthm0}
Let $T$ be an ordered tree bigraph, of radius $r$, and with $t\ge 2$ vertices. For all  $\tau>0$ there exists $\delta>0$ with the following property.
Let $G$ be a bigraph not containing $T$, with $|V_1(G)|, |V_2(G)|\ge n$, and with at most $(1-\tau)|V_1(G)|\cdot|V_2(G)|$ edges.
Choose $K$ such that $t^{K^r}=n$. 
Then
there are subsets $Z_i\subseteq V_i(G)$ for $i = 1,2$, with $|Z_1|,|Z_2|\ge \delta nt^{-5K^{r-1}}$, such that $Z_1,Z_2$ are anticomplete.
\end{thm}
It also implies:
\begin{thm}\label{symthm}
Let $T$ be an ordered tree bigraph, of radius $r$, and with $t\ge 2 $ vertices.
Let $G$ be an ordered bigraph not containing $T$ or its bicomplement, with $|V_1(G)|, |V_2(G)|\ge n$.
Choose $K$ such that $t^{K^r}=n$.
Then there is a pure pair $(Z_1,Z_2)$ with $Z_i\subseteq V_i(G)$ and $|Z_i|\ge (16t^2)^{-t}nt^{-5K^{r-1}}$ for $i = 1,2$.
\end{thm}

Similarly, the second result, \ref{linmainthm}, implies:
\begin{thm}\label{linmainthm0}
Let $T$ be an ordered tree bigraph. For all  $c,\tau>0$ there exists $\delta>0$ with the following property.
Let $G$ be a bigraph not containing $T$, with $|V_1(G)|, |V_2(G)|\ge n$, and with at most $(1-\tau)|V_1(G)|\cdot|V_2(G)|$ edges.
Then
there are subsets $Z_i\subseteq V_i(G)$ for $i = 1,2$, either with $|Z_1|\ge \delta n$ and $|Z_2|\ge \delta n^{1-c}$, or with
$|Z_1|\ge \delta n^{1-c}$ and $|Z_2|\ge \delta n$, such that $Z_1,Z_2$ are anticomplete.
\end{thm}
It also implies:
\begin{thm}\label{linsymthm}
Let $T$ be an ordered tree bigraph. For all  $c>0$ there exists $\delta>0$ with the following property.
Let $G$ be an ordered bigraph not containing $T$ or its bicomplement, with $|V_1(G)|, |V_2(G)|\ge n$.
Then there is a pure pair $(Z_1,Z_2)$ with $Z_i\subseteq V_i(G)$ for $i = 1,2$, either with $|Z_1| \ge \delta n$ and $|Z_2|\ge \delta n^{1-c}$, or with $|Z_1| \ge \delta n^{1-c}$ and $|Z_2|\ge \delta n$.
\end{thm}

The last six theorems all imply that $|Z_1|,|Z_2|=n^{1-o(1)}$.
It is easy to show, with a random graph argument, that this characterizes forests and their bicomplements. Indeed,
if $T$ is an (unordered) bigraph such that neither $T$ nor its bicomplement are forests,
then the conclusions of \ref{symthm} and \ref{linsymthm} (with ``ordered'' deleted)  are far from true; and therefore 
for ordered bigraphs they are at least as far from true.
More exactly, here is a standard example:
\begin{thm}\label{example}
Let $T$ be a bigraph, such that both $T$ and its bicomplement 
have a cycle of length at most $g$, and let $c>1-1/g$.
Then there is a bigraph $G$, with $|V_1(G)|, |V_2(G)|=n$,
which contains neither $T$ nor its bicomplement, and such that 
$\min(|Z_1|, |Z_2|)\le n^{c}$ for every pure
pair $(Z_1,Z_2)$ with $Z_i\subseteq V_i(G)$
for $i = 1,2$. 
\end{thm}
\Proof
Take $n$ large, and let $V_1,V_2$ be disjoint sets of cardinality $2n$; and for each $v_1\in V_1$ and $v_2\in V_2$, make $v_1,v_2$
adjacent independently, with probability $\frac12 n^{1/g-1}$. Then with high probability, there are fewer than $n/2$ cycles of length 
at most $g$, and no pure pair $Z_1,Z_2$ with $Z_i\subseteq V_i$ and $|Z_i|\ge n^c$ for $i = 1,2$. By deleting half of $V_1,V_2$
appropriately, 
we obtain a bigraph $G$ with girth more than $g$, which therefore does not contain $T$ or its bicomplement. This proves \ref{example}.~\bbox

\section{Reduction to the sparse case}
In this section we do two things. First, we deduce \ref{symthm} assuming \ref{mainthm}, and will prove the latter in the 
next section; and second, we deduce \ref{mainthm0} and \ref{linmainthm0} from \ref{mainthm} and \ref{linmainthm}.

We need the following lemma,
a version of a theorem of Erd\H{o}s, Hajnal and Pach~\cite{erdos} adapted for ordered bipartite graphs. (It is similar to a result of \cite{pure4} but with different parameters).
\begin{thm}\label{sparsify}
Let $H$ be an ordered bigraph, let $|V(H_i)|=h_i$ for $i = 1,2$, let $0<\vare<1/8$, let $d= \lceil 1/(4\vare)\rceil$, 
and let $m_1,m_2>0$ be integers. 
Let $G$ be an ordered bigraph not containing $H$, with $|V_1(G)|\ge h_1d^{h_2}m_1$ and $|V_2(G)|\ge 2h_1h_2m_2$. 
Then there are subsets $Y_i\subseteq V_i(G)$
with $|Y_i|=m_i$ for $i = 1,2$, such that either 
\begin{itemize}
\item every vertex in $Y_1$ has at most $\vare |Y_2|$ neighbours in $Y_2$, and every vertex in $Y_2$ has at most $\vare|Y_1|$ 
neighbours in $Y_1$, or
\item every vertex in $Y_1$ has at most $\vare |Y_2|$ non-neighbours in $Y_2$,
and every vertex in $Y_2$ has at most $\vare|Y_1|$ 
non-neighbours in $Y_1$.
\end{itemize}
\end{thm}
\Proof 
Divide $V_1(G)$ into $h_1$ disjoint intervals, each of 
cardinality at least $m_1d^{h_2}$, numbered
$B_u\;(u\in V_1(H))$ in order. Divide $V_2(G)$ into disjoint intervals $B_u\;(u\in V_2(H))$ of cardinality at least $m_2h_1$.
Choose $W\subseteq V_2(H)$
maximal such that for each $v\in W$ there exists $x_v\in B_v$, and for each $u\in V_1(H)$ there exist $Q_u\subseteq B_u$,
 with the following properties:
\begin{itemize}
\item $|Q_u|\ge m_1d^{h_2-|W|}$ for each $u\in V_1(H)$;
\item for each $u\in V_1(H)$ and each $v\in W$, if $u,v$ are $H$-adjacent then $x_v$ is complete to $Q_u$, and if $u,v$ are 
not $H$-adjacent then $x_v$ is anticomplete to $Q_u$.
\end{itemize}
This is possible since we may take $W=\emptyset$ and $Q_u=B_u$ for each $u\in V_1(H)$. Since $G$ does not contain $H$, it follows that $W\ne V_2(H)$. 
Choose $v\in V_2(H)\setminus W$.
Say $u\in V_1(H)$ is a {\em problem} for $x\in B_v$ if either $u,v$ are $H$-adjacent and 
$x$ has fewer than $|Q_u|/d$ neighbours in $Q_u$, or $u,v$ are not $H$-adjacent and $x$ has fewer than $|Q_u|/d$
non-neighbours in $Q_u$. From the maximality of $W$, for each $x\in B_v$ there exists $u\in V_1(H)$ that is a problem for $x$.
Since there are only $h_1$ possible problems, 
there exist $u\in V_1(H)$, and $C\subseteq B_v$ with $|C|\ge |B_v|/h_1$, such that 
for every $x\in C$, $u$ is a problem for $x$. By moving to the bicomplement if necessary, we may assume that $u,v$ are $H$-adjacent;
and so every vertex in $C$ has fewer than $|Q_u|/d$ neighbours in $Q_u$. Since $|Q_u|\ge m_1d^{h_2-|W|}\ge m_1d\ge 2m_1$ and 
$|C|\ge |B_v|/h_1\ge 2m_2$,
it follows by averaging that there are subsets $X_1\subseteq Q_u$ and $X_2\subseteq C$, of cardinality exactly $2m_1,2m_2$ respectively,
such that there are at most $|X_1|\cdot|X_2|/d=4m_1m_2/d$ edges joining them. Let $Y_1$ be the set of the 
$m_1$ vertices in $X_1$ that have fewest neighbours
in $X_2$; then they each have at most $4m_2/d$ neighbours in $X_2$. Define $Y_2$ similarly; then $|Y_i|=m_i$ for $i = 1,2$, and every
vertex in $Y_1$ has at most $4m_2/d\le \vare m_2$ neighbours in $Y_2$ and vice versa. This proves \ref{sparsify}.~\bbox

\noindent{\bf Proof of \ref{symthm}, assuming \ref{mainthm}.\ \ }
Let $T$ be an ordered tree bigraph, of radius $r$, and with $t$ vertices. 
Let $\vare = 1/(4t^2)$, and let $d=\lceil 1/(4\vare)\rceil$. Let $c=(16t^2)^{-t}$.
Let $G$ be a bigraph not containing $T$ or its bicomplement, with $|V_1(G)|, |V_2(G)|\ge n$. 
We may assume that $cn t^{-5K^{r-1}}>1$, for otherwise the result is true, taking $|Z_1|=|Z_2|=1$. 
Hence $2cn\ge 2t^{5K^{r-1}}\ge 1$.
Let $m$ be the largest integer
such that $m\le 2cn$. Thus $m\ge cn$, since $2cn\ge 1$.  

Let $|V_i(T)|=h_i$ for $i = 1,2$.
Now $n\ge  \max(h_1d^{h_2}m, 2h_1h_2m)$.
By \ref{sparsify}, and moving to the bicomplement if necessary, we may assume that there exist $Y_i\subseteq V_i(G)$
with $|Y_i|=m$ for $i = 1,2$, such that 
every vertex in $Y_1$ has at most $\vare |Y_2|$ neighbours in $Y_2$, and every vertex in $Y_2$ has at most $\vare|Y_1|$
neighbours in $Y_1$. Choose $J$ with $m=t^{J^r}$. By \ref{mainthm} applied to the ordered bigraph induced on $Y_1\cup Y_2$, there exist $Z_i\subseteq Y_i$
with $|Z_1|,|Z_2|\ge mt^{-5J^{r-1}}$, such that $Z_1,Z_2$ are anticomplete, where $m=t^{J^r}$. 
Since $J\le K$, it follows that  $|Z_1|,|Z_2|\ge cnt^{-5K^{r-1}}$. This proves \ref{symthm}.~\bbox

The proof that \ref{linmainthm} implies \ref{linsymthm} is similar and we omit it.

The result \ref{stars} of Kor\'andi, Pach, and Tomon~\cite{pach} has as a hypothesis that the bicomplement of $G$
has at least $\tau n^2$ edges. This is apparently much weaker than the hypothesis that every vertex of $G$
has degree at most $\vare n$, but in fact the ``not very dense'' hypothesis is as good as
the ``very sparse'' hypothesis, because of the next result, proved in \cite{pure4}.

\begin{thm}\label{countnonedges}
For all $c,\vare,\tau>0$ with $\vare<\tau$, there exists $\delta>0$ with the following property. Let $G$
be a bigraph with at most $(1-\tau)|V_1(G)|\cdot|V_2(G)|$  edges and with $V_1(G), V_2(G)\ne \emptyset$. Then
there exist $Z_i\subseteq V_i(G)$ with  $|Z_i|\ge \delta|V_i(G)|$ for $i = 1,2$, such that
there are fewer than
$(1-\vare)|Y_1|\cdot|Y_2|$ edges between $Y_1,Y_2$ for all subsets $Y_i\subseteq Z_i$ with
$|Y_i|\ge  c|Z_i|$ for $i = 1,2$.
\end{thm}

We deduce:

\begin{thm}\label{countnonedges2}
For every ordered bigraph $H$, and for all $\vare,\tau>0$, there exists $\delta>0$ with the following property. Let $G$
be an ordered bigraph not containing $H$, with at most $(1-\tau)|V_1(G)|\cdot|V_2(G)|$  edges. Then there exist
$Y_i\subseteq V_i(G)$ with $|Y_i|\ge \delta|V_i(G)|$ for $i = 1,2$, such that every vertex in $Y_1$ has at most $\vare |Y_2|$ 
neighbours in $Y_2$, and every vertex in $Y_2$ has at most $\vare|Y_1|$
neighbours in $Y_1$.
\end{thm}
\Proof We may assume that $\vare<1/8$ and $\vare<\tau$, by reducing $\vare$, and $|V_1(H)|, |V_2(H)|\ne \emptyset$, by adding
vertices to $H$. Let
$h_i=|V(H_i)|$ for $i = 1,2$, let 
let $d= \lceil 1/(4\vare)\rceil$, and let $1/c=\max(h_1d^{h_2}, 2h_1h_2)$.
Choose $\delta'$ such that \ref{countnonedges} holds with $\delta$ replaced by $\delta'$, and let $\delta=c\delta'$.
Now let $G$
be an ordered bigraph not containing $H$, with at most $(1-\tau)|V_1(G)|\cdot|V_2(G)|$  edges. 
We may assume that $V_1(G), V_2(G)\ne \emptyset$. By \ref{countnonedges},
there exist $Z_i\subseteq V_i(G)$ with  $|Z_i|\ge \delta'|V_i(G)|$ for $i = 1,2$, such that
there are fewer than
$(1-\vare)|Y_1|\cdot|Y_2|$ edges between $Y_1,Y_2$ for all subsets $Y_i\subseteq Z_i$ with
$|Y_i|\ge  c|Z_i|$ for $i = 1,2$. By \ref{sparsify}, applied to the ordered sub-bigraph induced on $Z_1\cup Z_2$,
there exist $Y_i\subseteq Z_i$ with $|Y_i|\ge c|V_i(G)|$ for $i = 1,2$, such that either
\begin{itemize}
\item every vertex in $Y_1$ has at most $\vare |Y_2|$ neighbours in $Y_2$, and every vertex in $Y_2$ has at most $\vare|Y_1|$
neighbours in $Y_1$, or
\item every vertex in $Y_1$ has at most $\vare |Y_2|$ non-neighbours in $Y_2$,
and every vertex in $Y_2$ has at most $\vare|Y_1|$
non-neighbours in $Y_1$.
\end{itemize}
Since there are fewer than $(1-\vare)|Y_1|\cdot|Y_2|$ edges between $Y_1,Y_2$, the second is impossible, and so the 
first holds. Then for $i = 1,2$, $|Y_i|\ge c|Z_i|\ge c\delta'|V_i(G)|=\delta|V_i(G)|$. This proves \ref{countnonedges2}.~\bbox

\noindent{\bf Proof of \ref{mainthm0}, assuming \ref{mainthm}.\ \ }
Let $T$ be an ordered tree bigraph, of radius $r$, and with $t\ge 2$ vertices, and let $\tau>0$.
Let $\vare=1/(4t^2)$, and choose $\delta>0$ as in \ref{countnonedges}, with $H$ replaced by $T$.

Let $G$ be a bigraph not containing $T$, with $|V_1(G)|, |V_2(G)|\ge n$, such that the bicomplement of $G$ has at least
$\tau|V_1(G)|\cdot|V_2(G)|$ edges. From the choice of $\delta$, there exist
$Y_i\subseteq V_i(G)$ with $|Y_i|\ge \delta|V_i(G)|$ for $i = 1,2$, such that every vertex in $Y_1$ has at most $\vare |Y_2|$
neighbours in $Y_2$, and every vertex in $Y_2$ has at most $\vare|Y_1|$
neighbours in $Y_1$. By \ref{mainthm} applied to the sub-bigraph $G[Y_1\cup Y_2]$, 
there are subsets $Z_i\subseteq Y_i$ for $i = 1,2$, with $|Z_1|,|Z_2|\ge \delta nt^{-5k^{r-1}}$, such that $Z_1,Z_2$ are anticomplete,
where $k$ satisfies $t^{k^r}=\delta n$. Choose $K$ such that $t^{K^r}=n$; thus $K\ge k$, and so $|Z_1|,|Z_2|\ge \delta nt^{-5K^{r-1}}$.
This proves \ref{mainthm0}.~\bbox

The proof that \ref{linmainthm} implies \ref{linmainthm0} is similar and we omit it.

\section{Proof of the first main theorem}
In this section we prove \ref{mainthm}, which we restate:
\begin{thm}\label{mainthm2}
Let $T$ be a ordered tree bigraph, of radius $r$, and with $t\ge 2$ vertices. 
Let $G$ be a bigraph  with $|V_1(G)|, |V_2(G)|\ge n$, that does not contain $T$, and such that every vertex has degree
at most $n/(4t^2)$. Choose $K$ such that $t^{K^r}=n$. Then
there are two anticomplete subsets $Z_i\subseteq V_i(G)$ for $i = 1,2$, with $|Z_1|,|Z_2|\ge nt^{-5K^{r-1}}$.
\end{thm}

\Proof 
Let $\vare=1/(4t^2)$.
Since $\vare<1$, $G$ is not complete bipartite, and so we may assume that $t^{-5K^{r-1}}n>1$ (or else the theorem holds);
that is, $t^{5K^{r-1}}<t^{K^r}$. So $K\ge 5$, and $n\ge t^{5^r}$.

If $r=1$, then $T$ has a vertex of degree $d$ say, and all other vertices of $T$ are neighbours of $d$. Let this vertex belong to 
$V_1(T)$ say. Since $G$ does not contain $T$, all vertices in $V_1(G)$ have degree at most $d-1$. Choose a set $Z_1$ of at most
$n/d$ vertices in $V_1(G)$; then the set of vertices with neighbours in $X$ has cardinality at most $(d-1)|X|\le (d-1)n/d$,
and so there is a set of at least $n/d$ vertices in $V_2(G)$ anticomplete to $Z_1$. So to prove \ref{mainthm} in this case,
we just have to check that $\lfloor n/d \rfloor\ge nt^{-5K^{r-1}}$. But $t=d+1$ and $n\ge t$ (because $n\ge t^{5^r}$), so
$\lfloor n/d \rfloor\ge n/(2(t-1))$; and now it remains to check that $n/(2(t-1))\ge nt^{-5K^{r-1}}$, which is clear. Thus we may assume that $r\ge 2$, and so $t\ge 4$.

Choose a real number $x\ge 0$ with $x\le K^{r-1}$, 
maximum such that there exist $A_1\subseteq V_1(G)$ and $A_2\subseteq V_2(G)$
with the properties that
\begin{itemize}
\item $|A_1|,|A_2|\ge nt^{-x}$;
\item every vertex in $A_1$ has at most $\vare nt^{-Kx}$ neighbours in $A_2$ and vice versa.
\end{itemize}
This is possible since we may take $x=0$ and $A_i=V_i(G)$ for $i = 1,2$.
Let $A_1,A_2$ be as above. Let $d=\vare nt^{-Kx}$.
Since $|A_1|\ge nt^{-x}\ge nt^{-5K^{r-1}}$, we may assume that $A_1$ is not anticomplete to $A_2$ (or else the theorem holds); so $d\ge 1$.
For $1\le s\le r-1$, let $k_s= 4(K^{s-1}+K^{s-2}+\cdots + 1)$.
\\
\\
(1) {\em $x< K^{r-1}-k_{r-1}$.}
\\
\\
Suppose not.
Since $|A_1|\ge t^{-x}n\ge t^{-K^{r-1}}n\ge 2t^{-5K^{r-1}} n$, there exists a set $X\subseteq A_1$ of cardinality $\lceil t^{-5K^{r-1}} n\rceil\le 2t^{-5K^{r-1}} n$.
The union of the neighbours in $A_2$ of vertices in $X$ has cardinality at most 
$(2t^{-5K^{r-1}} n)(\vare nt^{-Kx})$.
Since $|A_2|/2\ge t^{-5K^{r-1}} n$, there are fewer than $|A_2|/2$ vertices in $A_2$ anticomplete to $X$; 
so $(2t^{-5K^{r-1}} n)(\vare nt^{-Kx})\ge  |A_2|/2$, and hence
$4\vare t^{-5K^{r-1}} n t^{-Kx}\ge t^{-x}$. Consequently
$4\vare t^{-5K^{r-1}} n \ge t^{(K-1)x}$.
But 
$$(K-1)x\ge (K-1) (K^{r-1}-k_{r-1})=K^r-5K^{r-1}+4,$$
and so 
$$4\vare t^{-5K^{r-1}} n \ge t^{K^{r}-5K^{r-1}+4}.$$
Hence 
$n \ge t^{K^{r}+4},$ a contradiction.  
This proves (1).

\bigskip
Since $x<K^{r-1}-4$, and $n>t^{5K^{r-1}}$, it follows that $nt^{-x-1}\ge t$. Hence 
$$|A_i|\ge nt^{-x}\ge (t-1) nt^{-x-1} +t$$
for $i = 1,2$.
Since $|V(T_1)|\le t-1$, it follows that we may choose $|V_1(T)|$ disjoint blocks $B_u\;(u\in V_1(T))$, 
all intervals of $A_1$, and of 
cardinality $\lceil nt^{-x-1}\rceil$, numbered in order. Partition $A_2$ into blocks $B_v\;(v\in V_2(T))$ similarly.
For $1\le s\le r-1$, let
$p_s = dt^{-Kk_s}$. Let $p_r=1$.  For $2\le s\le r$, let $f_s=dt^2/p_{s-1}=t^{Kk_{s-1}+2}$. So $f_s\ge t^{4K+2}$.

Since $T$ has radius at most $r$, there is a vertex $v_0\in V(T)$ such that every vertex of $T$ can be joined to $v_0$
by a path of length at most $r$.
For $1\le s\le r$ let $T_s$ be the subtree of $T$ induced on the vertices with distance at
most $s$ from $t_0$. So $V(T_0) = \{t_0\}$, and $T_r=T$. Let $L_s$ be the set of vertices with distance exactly $s$ from $v_0$.
For $2\le s\le r$, and  each edge $uv$ of $T$ with $u\in L_{s-1}$
and $v\in L_s$, choose $X_{uv}\subseteq B_u$
and $Y_{uv}\subseteq B_v$ satisfying the following conditions:
\begin{itemize}
\item every vertex in $X_{uv}$ has fewer than $p_s$ neighbours in $B_v\setminus Y_{uv}$;
\item $|Y_{uv}|\le f_s|X_{uv}|$ and $|Y_{uv}|\le |B_v|/2$; and
\item subject to these conditions, $Y_{uv}$ is maximal.
\end{itemize}
This is possible since we could take $X(uv)=Y(uv)=\emptyset$ to satify the first two bullets.
\\
\\
(2) {\em For $2\le s\le r$, and  each edge $uv$ of $T$ with $u\in L_{s-1}$
and $v\in L_s$, we may assume that $|X_{uv}|=\lceil |Y_{uv}|/f_s\rceil$, and $|X_{uv}|\le |B_u|/(2t)$, and $|Y_{uv}|\le |B_v|/2-td$.}
\\
\\
We may assume that $|X_{uv}|=\lceil |Y_{uv}|/f_s\rceil$, by removing elements from $X_{uv}$ if necessary. Suppose first that $s=r$. 
Then $X_{uv}$
is anticomplete to $B_v\setminus Y_{uv}$ (because $p_r=1$), and so either $|X_{uv}|<t^{-5K^{r-1}} n$ or 
$|B_v\setminus Y_{uv}|<t^{-5K^{r-1}} n$. 
The second implies
that $|B_v|<2t^{-5K^{r-1}} n$ (since $|Y_{uv}|\le |B_v|/2$), and so $nt^{-x-1}<2t^{-5K^{r-1}} n$, that is,
$t^{5K^{r-1}-x-1}<2$, a contradiction. So $|X_{uv}|<t^{-5K^{r-1}} n$. Since $t^{-5K^{r-1}} n\le nt^{-x-1}/(2t)$, it follows that
$|X_{uv}|\le |B_u|/(2t)$.
Also, 
$|Y_{uv}|\le f_rt^{-5K^{r-1}} n$. We claim that $f_rt^{-5K^{r-1}} n\le |B_v|/2-td$. Suppose not; then either $f_rt^{-5K^{r-1}}n>|B_v|/4$ or $td>|B_v|/4$.
The first implies that $t^{Kk_{r-1}+2}t^{-5K^{r-1}} n>nt^{-x-1}/4$, and so $t^{Kk_{r-1}+x+4}t^{-5K^{r-1}} >1$.
Hence $Kk_{r-1}+x+4-5K^{r-1}>0$. But $x\le K^{r-1}-k_{r-1}$ by (1), so
$$Kk_{r-1}+ K^{r-1}-k_{r-1}+ 4-5K^{r-1}>0,$$
a contradiction (in fact, the left side sums to zero).
The second implies
that $4\vare t^{-Kx}>t^{-x-2}$ and so $K< 1$ since $4\vare =t^{-2}$, a contradiction. Thus when $s=r$, all three statements of (2) hold.

Now we assume that $2\le s<r$. We have $|X_{uv}|\le |Y_{uv}|/f_s+1\le |B_{v}|/(2f_s) +1\le |B_v|/4$, because
$f_s\ge t^{4K+2}\ge 4$ and $|B_v|\ge 8$. There are at most $|X_{uv}|p_s$ edges between $X_{uv}$ and $B_v\setminus Y_{uv}$;
and so at most $|X_{uv}|$ vertices in $B_v\setminus Y_{uv}$ have at least $p_s$ neighbours in $X_{uv}$. Since
$|B_v\setminus Y_{uv}|\ge |B_v|/2\ge 2|X_{uv}|$, 
the $|X_{uv}|$ vertices in $B_v\setminus Y_{uv}$ with fewest neighbours in $X_{uv}$
each have fewer than $p_s$ neighbours in $X_{uv}$. Since $p_s=dt^{-Kk_s}$, and $k_s\le 4K^{r-2}+4K^{r-3}+\cdots+4$, the maximality of $x$ 
and (1) imply that $|X_{uv}|<nt^{-x-k_s}$. Hence $|X_{uv}|\le |B_u|/(2t)$, because $nt^{-x-k_s}\le nt^{-x-1}/(2t)$
(because $k_s\ge 3$). Thus the second claim holds. For the third claim, since $|Y_{uv}|\le f_s|X_{uv}|$, it suffices to show that
$f_s|X_{uv}|\le |B_v|/2-td$, and to prove this, it suffices to show that $f_s|X_{uv}|\le |B_v|/4$ and $td\le  |B_v|/4$.
To show the first, it suffices to show that $nt^{-x-k_s}f_s\le |B_v|/t$, that is,
$nt^{-x-k_s}t^{Kk_{s-1}+2}\le nt^{-x-2}$,
which simplifies to
$4+ Kk_{s-1}\le k_s$, and this holds with equality. To show that $td\le  |B_v|/4$, it suffices to show that
$t\vare n t^{-Kx}\le nt^{-x-1}/4$, which simplifies to $4\vare t^{2+x-Kx}\le 1$; and this is true since $4\vare t^2\le 1$, 
and $K\ge 1$. This proves (2).

\bigskip

For $2\le s\le r$, and each $u\in L_{s-1}$, let $X_u$ be the union of the sets $X_{uv}$ over all $v\in L_{s}$ that are $T$-adjacent to $u$. 
Then:
\\
\\
(3) {\em For $2\le s\le r$, and each $u\in L_{s-1}$, $|X_u|\le |B_u|/2$.}
\\
\\
For each $v\in L_{s}$ that is $T$-adjacent to $u$, 
$|X_{uv}|\le |B_u|/(2t)$ by (2), and the claim follows. This proves (3).

\bigskip

Let $P_{v_0}=B_{v_0}$. For $s=1\ll r-1$ we will choose $P_v\subseteq B_v\setminus X_v$ for each $v\in L_s$, and 
$y_v\in P_v$ for each $v\in L_{s-1}$, satisfying the following conditions:
\begin{itemize}
\item for all distinct $u,v\in V(T_{s-1})$, $u,v$ are $T$-adjacent if and only if $y_u, y_v$ are $G$-adjacent;
\item for all $u\in V(T_{s-1})$ and $v\in L_s$, and all $y\in P_v$, $u,v$ are $T$-adjacent if and only if $y_u, y$ are $G$-adjacent;
\item for each $v\in L_s$, $|P_v|\ge p_s$.
\end{itemize}
First let us assume $s=1<r$. If there exists $y\in B_{v_0}$ with at least $p_1$ neighbours in $B_v\setminus X_v$ for each $v$ 
such that $v_0, v$ are $T$-adjacent,
then we may set $y_{v_0}=y$; so we assume there is no such $y$. Consequently there is a $T$-neighbour $v$ of $v_0$ such that for at 
least $|B_{v_0}|/t\ge nt^{-x-2}$ vertices $y\in B_{v_0}$, $y$ has fewer than $p_1$ neighbours in $B_v\setminus X_v$. Choose a set $X$ of 
exactly $\lceil nt^{-x-2} \rceil$ such vertices $y$. Since 
$|X|\le nt^{-x-2}+1\le 2nt^{-x-2}$,  and $|B_{v_0}|\ge 8nt^{-x-2}$ (because $t\ge 4$), it follows that 
$|X|\le |B_{v_0}|/4$. Since $|X_v|\le |B_v|/2$, it follows that $|B_v\setminus X_v|\ge 2|X|$, and so at least $|X|$ vertices in 
$B_v\setminus X_v$ have at most $p_1$ neighbours in $X$. Since $p_1=dt^{-Kk_1}=dt^{-4K}$, the maximality of $x$
implies that $|X|<nt^{-x-4}$, a contradiction.
So we can satisfy the three bullets above when $s=1\le r-1$. 

Suppose that $2\le s\le r$, and
we have chosen
$P_v\subseteq B_v\setminus X_v$ for each $v\in L_{s-1}$ and $y_v\in P_v$ for each $v\in V(T_{s-2})$. We must define 
$P_v\subseteq B_v\setminus X_v$ for each $v\in L_{s}$, and
$y_v\in P_v$ for each $v\in L_{s-1}$, satisfying the bullets above. From the symmetry we may assume that
$L_s\subseteq V_1(T)$.

Let $C$ be the set of vertices in $A_2$ that are equal or adjacent to $y_v$ for some $v\in V(T_{s-2})$. 
Let $y_u\in P_u$ for each $u\in L_{s-1}$; we call $(y_u:u\in L_{s-1})$ a {\em transversal}. A transversal
$(y_u:u\in L_{s-1})$ is {\em valid} if 
for each edge $uv$ of $T$ with $u\in L_{s-1}$ and $v\in L_s$, there are at least $p_s$ vertices in $B_v\setminus C$
that are adjacent to $y_u$ and that have no other neighbour in $\{y_{u'}:\;u'\in L_{s-1}\}$.
\\
\\
(4) {\em There is a valid transversal.}
\\
\\
Suppose not. Let $E$ be the set of ordered pairs $(u,v)$ such that $uv$ is an edge of $T$ with $u\in L_{s-1}$ and $v\in L_s$.
Then for every transversal $(y_u:u\in L_{s-1})$, there exists
$(u,v)\in E$ such that there are fewer than $p_s$ vertices in $B_v\setminus C$
that are adjacent to $y_u$ and that have no other neighbour in $\{y_{u'}:\;u'\in L_{s-1}\}$. Call $(u,v)$ a {\em problem}
for the transversal $(y_u:u\in L_{s-1})$. Since $|E|=|L_s|$, there are only $|L_s|$ possible problems, and so there 
exists $(u,v)\in E$ that is a problem for at least a fraction $1/|L_s|$ of all transversals. Hence there exist
a subset  $X\subseteq P_u$ with $|X|\ge |P_u|/|L_s|\ge p_{s-1}/t$, and a choice of $y_{u'}\in P_{u'}$
for each $u'\in L_{s-1}\setminus \{u\}$, 
such that for all
$y_u\in X$, $(u,v)$ is a problem for the transversal $(y_{u'}:u'\in L_{s-1})$.
Let $C'$ be the set of vertices in $A_2$ that are adjacent to
a vertex in $(y_{u'}:u'\in L_{s-1}\setminus \{u\})$. Since every vertex in $C\cup C'$ has a neighbour $y_w$
for some $w\in V(T_{s-1})\setminus \{u\}$, it follows that $|C\cup C'|\le dt$. Every vertex in $X$ has fewer than $p_s$
neighbours in $B_v\setminus (C\cup C')$. If $(C\cup C')\cap B_v\not\subseteq Y_{uv}$
let $Y=(C\cup C')\cap B_v$, and otherwise let $Y$ be a singleton subset of $B_v\setminus Y_{uv}$.
Thus $|Y|\le dt\le f_s|X|$, since $dt\ge 1$ by (1). Consequently 
$|Y\cup Y_{uv}|\le f_s|X\cup X_{uv}|$, since $X\cap X_{uv}=\emptyset$; and since $|C\cup C'\cup Y_{uv}|\le |B_v|/2$
(because $|C\cup C'|\le dt$ and by (2)), and every vertex in $X\cup X_{uv}$ has fewer than $p_s$
neighbours in $B_v\setminus (Y\cup Y_{uv})$, this contradicts the maximality of $Y_{uv}$. This proves (4).

\bigskip
From (4), the inductive definition of $y_v\;(v \in V(T_{r-1}))$ and $P_v\;(v\in V(T))$ is complete.
For each $v\in L_s$, choose $y_v\in P_v$. Then the map sending each $v\in V(T)$ to $y_v$ is 
an ordered parity-preserving isomorphism of $T$ to an induced subgraph of $G$,
a contradiction. This proves \ref{mainthm2}.~\bbox

\section{Parades}

Now we begin the proof of \ref{linmainthm}, the second main result mentioned in the introduction. 
This proof was derived from, and still has some ingredients in common with, the proof of the main theorem of \cite{pure6}, but 
it has needed some serious modification, in order to persuade one of the two sets $Z_1,Z_2$ to be linear.

Let $G$ be a bigraph (not necessarily ordered), and let $I$ be a set of nonzero integers. We denote 
$I^+=\{i\in I: i>0\}$ and $I^-=I\setminus I^+$.
Let the sets $B_i\;(i\in I)$ be nonempty, pairwise disjoint subsets of $V(G)$, such that $B_i\subseteq V_1(G)$ if $i<0$
and $B_i\subseteq V_2(G)$ if $i>0$.
We call 
$\mathcal{P}=(B_i:i\in I)$ a {\em parade} in $G$. 
Its {\em length} is the 
pair $(|I^-|,|I^+|)$, and its {\em width}
is the pair $(w_1,w_2)$ where $w_1= \min(|B_i|:i\in I^-)$ and $w_2= \min(|B_i|:i\in I^+)$, taking $w_k=|V_k(G)|$ if the corresponding set $I^-$ or $I^+$ is empty. 
We call the sets $B_i$ 
the {\em blocks} of the parade. What matters is that the blocks are not too small. 
(We used 
the same word in~\cite{pure4} for a similar but slightly different object.)

If $I'\subseteq I$, then $(B_i:i\in I')$ is a parade, called a {\em sub-parade}
of $\mathcal{P}$. If $B_i'\subseteq B_i$ is nonempty for each $i\in I$, 
then $(B_i':i\in I)$ is a parade, called a {\em contraction} of $\mathcal{P}$. 

Let $X,Y$ be disjoint nonempty subsets of $V(G)$. The {\em max-degree from $X$ to $Y$} is defined to be the maximum 
over all $v\in X$ of the number of neighbours of $v$ in $Y$.
Let $(B_i:i\in I)$ be a parade in a bigraph $G$, and for all $i,j\in I$ of opposite sign, 
let $d_{i,j}$ be the max-degree from $B_i$ to $B_j$. (For all other pairs $i,j$ we define $d_{i,j}=0$.)
We call $d_{i,j}\;(i,j\in I)$ the {\em max-degree function} of the parade.
The product of the numbers $d_{j,h}$ for all pairs $h,j$ where $h\in I^-$ and $j\in I^+$ is called the {\em max-degree product} of $\mathcal{B}$.
We just need this ``product'' definition for the next theorem.

Let $0<\phi,\mu$. We say that $\mathcal{B}$ is
{\em $(\phi,\mu)$-shrink-resistant} if for all $h\in I^-$ and $j\in I^+$, and 
for all $X\subseteq B_h$ and $Y\subseteq B_j$ with $|X|\ge \mu|B_h|$ and $|Y|\ge \mu|B_j|$,
the max-degree from $Y$ to $X$ is more than 
$d_{j,h}|V_1(G)|^{-\phi}$.
We begin with:

\begin{thm}\label{shrinkresistant}
Let $\mathcal{B}=(B_i:i\in I)$ be a
parade in a bigraph $G$, and 
let $0<  \phi,\mu$ with $\mu\le 1$. Let $\beta=\mu^{1+|I|^2/\phi}$.   Then either
\begin{itemize}
\item there exist $h\in I^-$ and $j\in I^+$, and $X\subseteq B_h$ and $Y\subseteq B_j$ with 
$\frac{|X|}{|B_h|}, \frac{|Y|}{|B_j|}\ge \beta$, such that $X, Y$ are anticomplete; or
\item there is a  $(\phi,\mu)$-shrink-resistant contraction $(B_i':i\in I)$ of $\mathcal{B}$, such that 
$|B_i'|\ge \beta |B_i|$
for each $i\in I$.
\end{itemize}
\end{thm}
Let $S=\lfloor |I|^2/\phi \rfloor$.
Choose an integer $s$ with $0\le s\le S+1$ and with $s$ maximum such that there is
a contraction $\mathcal{B}'=(B_i':i\in I)$ of $\mathcal{B}$ with
\begin{itemize}
\item $|B_i'|\ge \mu^s|B_i|$ for each $i\in I^-$; and
\item max-degree product at most $|V_1(G)|^{|I|^2-\phi s}$.
\end{itemize}
(This is possible since we may take $s=0$ and
$\mathcal{B}'=\mathcal{B}$.) Let $d_{h,j}\;(h,j\in I)$ be the max-degree function of $\mathcal{B}'$.
\\
\\
(1) {\em We may assume that $d_{j,h}\ge 1$ for all $h\in I^-$ and $j\in I^+$, and so $s\le S$.}
\\
\\
If $d_{j,h}<1$, then $d_{j,h}=0$, since it is
an integer. Thus $B_h', B_j'$ are anticomplete.
Since $s\le S+1$ and hence $\mu^s\ge \mu^{S+1}\ge \beta$, it follows that $|B_h'|/|B_h|, |B_j'|/|B_j|\ge \beta$, and the
first outcome of the theorem holds. Thus we may assume that $d_{j,h}\ge 1$. Hence the
max-degree product of $\mathcal{B}'$ is at least one, and since it is at most
$|V_1(G)|^{|I|^2-\phi s}$, it follows that $|I|^2-\phi s\ge 0$.
Hence $s\le S$. This proves (1).
\\
\\
(2) {\em $(B_i':i\in I)$ is $(\phi,\mu)$-shrink-resistant.}
\\
\\
Let $h\in I^-$ and $j\in I^+$, and let $C_h\subseteq B_h'$ and $C_j\subseteq B_j'$, with
$|C_h|\ge \mu|B_h'|$ and $|C_j|\ge \mu |B_j'|$. For all $i\in I$
with $i\ne h,j$ let $C_i=B_i'$, and let $d$ be the max-degree from $C_j$ to $C_h$.
From the maximality of $s$, and since $s\le S$, it follows that the
max-degree product of $(C_i:i\in I)$ is more than $|V_1(G)|^{|I|^2-\phi (s+1)}$. 
Since the first is at most $d/d_{h,j}$ times the
max-degree product of $(B_i':i\in I)$, which is at most $|V_1(G)|^{|I|^2-\phi s}$, it follows that
$d/d_{h,j}> |V_1(G)|^{-\phi}$.
This proves (2).

\bigskip

Since 
$|B_i'|\ge \mu^S|B_i|\ge \beta|B_i|$ for each $i\in I$, the second outcome of the theorem holds. This proves
\ref{shrinkresistant}.~\bbox

Let $(B_i:i\in I)$ be a parade in a bigraph $G$, 
and let $0< \tau,\phi,\mu$. 
We say that $\tau$ is a 
{\em $(\phi,\mu)$-band} for $(B_i:i\in I)$ if for all $h\in I^-$ and $j\in I^+$:
\begin{itemize}
\item the max-degree from 
$B_j$ to $B_h$ is at most $\tau|B_h|$; and
\item for all $X\subseteq B_h$ and $Y\subseteq B_j$ with $|X|\ge \mu|B_h|$ and $|Y|\ge \mu|B_j|$,
the max-degree from $Y$ to $X$ is more 
than $\tau|V_1(G)|^{-\phi}|B_h|$.
\end{itemize}

\begin{thm}\label{ramsey}
Let $k\ge 0$ be an integer, and let $0< \phi, \mu$. Then there exists an integer $K\ge k$
with the following property. Let $G$ be a bigraph, 
and let $(B_i:i\in I)$ be a $(\phi,\mu)$-shrink-resistant parade in $G$, of length at least $(K,K)$. 
Then there exists $J\subseteq I$ with $|J^-|=|J^+|=k$ such that $(B_i:i\in J)$ 
has a $(2\phi,\mu)$-band.
\end{thm}
\Proof
Let $K\ge 1$ be an integer such that for every complete bipartite graph with bipartition $(H,J)$ where $|H|,|J|\ge K$,
and every colouring of its edges with $ \lfloor 1/\phi+1\rfloor$ colours, there exist $H'\subseteq H$ and $J'\subseteq J$
with $|H'|,|J'|= k$ such that all edges between $H', J'$ have the same colour.

Now let $(B_i:i\in I)$ be a $(\phi,\mu)$-shrink-resistant parade in $G$,
with max-degree function $d_{i,j}\;(i,j\in I)$. 

For all $h\in I^-$ and $j\in I^+$, there is an integer $s$ such that 
$$|V_1(G)|^{-(s+1)\phi}< d_{j,h}/|B_h|\le  |V_1(G)|^{-s\phi}.$$
We call $s$ the {\em type} of the pair $(h,j)$.
Since $|V_1(G)|^{-(s+1)\phi}< d_{j,h}/|B_h|\le 1$, it follows that $-(s+1)\phi< 0$, and so $s\ge 0$; and since
$1/|V_1(G)|\le d_{j,h}/|B_h|\le |V_1(G)|^{-s\phi}$ (because $d_{j,h}>0$ from the definition of $(\phi,\mu)$-shrink-resistant), 
it follows that $1\le |V_1(G)|^{1-s\phi}$, and so $s\le 1/\phi$.
Hence $s$ is one of the integers $0,1\ll \lfloor 1/\phi \rfloor$.
From the choice of $K$,
there exists $J\subseteq I$ with $|J^-|=|J^+|=k$, 
such that every pair $(h,j)$ with $h\in J^-$ and $j\in J^+$ has the same type, $s$ say. 
Let $\tau=|V_1(G)|^{-s\phi}$; then for all
$h\in J^-$ and $j\in J^+$,
$$\tau|V_1(G)|^{-\phi}< d_{j,h}/|B_h|\le \tau.$$
We claim that $\tau$ is a $(2\phi,\mu)$-band for $(B_i:i\in J)$.
To show this, it remains to show that 
for all $h\in J^-$ and $j\in J^+$, and
for all $X\subseteq B_h$ and $Y\subseteq B_j$ with $|X|\ge \mu|B_h|$ and $|Y|\ge \mu|B_j|$,
the max-degree from $Y$ to $X$ is more than $\tau|V_1(G)|^{-2\phi}|B_h|$.
But $\mathcal{B}$ is $(\phi,\mu)$-shrink-resistant, and so the max-degree from $Y$ to $X$ is more than $d_{j,h}|V_1(G)|^{-\phi}$;
and since $d_{j,h}\ge \tau|V_1(G)|^{-\phi}|B_h|$, the claim follows.
This proves \ref{ramsey}.~\bbox

By combining \ref{shrinkresistant} and \ref{ramsey}, we deduce:

\begin{thm}\label{homog}
Let $k\ge 0$ be an integer, and let  $0<  \phi,\mu$ with $\mu\le 1$. Then there exists an integer $K>0$ with the following property.
Let $\mathcal{B}=(B_i:i\in I)$ be a
parade  of length at least $(K,K)$ in a bigraph $G$.
Let $\beta=\mu^{1+2K^2/\phi}$.   Then either
\begin{itemize}
\item there exist $h\in I^-$ and $j\in I^+$, and $X\subseteq B_h$ and $Y\subseteq B_j$ with
$\frac{|X|}{|B_h|}, \frac{|Y|}{|B_j|}\ge \beta$, such that $X,Y$ are anticomplete; or
\item there exist $J\subseteq I$ with $|J^-|=|J^+|=k$, and a subset $B_i'\subseteq B_i$ with $|B_i'|\ge \beta |B_i|$ for each $i\in J$, 
such that 
$(B_i':i\in J)$ has a $(\phi,\mu)$-band.
\end{itemize}
\end{thm}
\Proof
Let $K$ satisfy \ref{ramsey} with $\phi$ replaced by $\phi/2$.
Let $G$ be a bigraph, and let $\mathcal{B}=(B_i:i\in I)$ be a
parade in $G$,  of length at least $(K,K)$. By \ref{shrinkresistant}, either
\begin{itemize}
\item there exist $h\in I^-$ and $j\in I^+$, and $X\subseteq B_i$ and $Y\subseteq B_j$ with
$|X|/|B_h|, |Y|/|B_j|\ge \beta$, such that $X, Y$ are anticomplete; or
\item there is a  $(\phi/2,\mu)$-shrink-resistant contraction $\mathcal{B}'=(B_i':i\in I)$ of $\mathcal{B}$, such that
$|B_i'|\ge \beta |B_i|$
for each $i\in I$.
\end{itemize}
In the first case the first outcome of the theorem holds. In the second case,
by \ref{ramsey} applied to $\mathcal{B}'$, the second outcome of the theorem holds.  This proves \ref{homog}.~\bbox

\section{Covering with leaves}

Again, this section concerns graphs rather than ordered graphs. If $G$ is a graph and 
$A,B\subseteq V(G)$ are disjoint, we say $A$ {\em covers} $B$ if every vertex of $B$
has a neighbour in $A$.

\begin{thm}\label{newleafcover}
Let $k\ge 1$ be an integer, and let $0< \tau, \phi,\mu$, with $\mu\le 1/(8k)$ and $\tau\le 1/(8k^2)$.
Let $G$ be a bigraph and let
let $\mathcal{A}=(A_i:i\in I)$ be a parade in $G$, with $|I^+|, |I^-|\le k$, such that
$\tau$ is a $(\phi,\mu)$-band for $\mathcal{A}$.
Then for each $h\in I^-$ there exist $B_h, C_h\subseteq A_h$, and for each $h\in I^-$ and $j\in I^+$ there 
exists $D_{h,j}\subseteq A_j$, with the following properties:
\begin{itemize}
\item 
$C_h\subseteq B_h$, and $|B_h|\ge |A_h|/2$, and 
$|C_h|\ge |V_1(G)|^{-k\phi}|A_h|/16$, for each $h\in I^-$; and
\item 
$D_{h,j}$ is anticomplete to $B_i$
for all $i\in I^-\setminus \{h\}$, and is 
anticomplete to $B_h\setminus C_h$, and covers $C_h$,
for each $h\in I^-$ and $j\in I^+$.
\end{itemize}
\end{thm}
\Proof
For each $h\in I^-$ and $j\in I^+$, every vertex in $A_j$ has at most $\tau|A_h|$ neighbours in $A_h$, and so there
are at most $\tau|A_h|\cdot|A_j|$ edges between $A_h, A_j$. Hence at most $1/(2k)$ vertices in $A_h$ have at least $2k\tau|A_j|$
neighbours in $A_j$. For each $h\in I^-$, let $P_h$ be the set of vertices $v\in A_h$ such that for each $j\in I^+$,
$v$ has fewer than 
$2k\tau|A_j|$ neighbours in $A_j$. It follows that $|P_h|\ge |A_h|/2$ for all $h\in I^-$.

Choose $H\subseteq I^-$ maximal such that for each $h\in H$ there exists $Q_h\subseteq P_h$ with 
$|Q_h|\ge |V_1(G)|^{-k\phi}|A_h|/8$, and for all $h\in H$ and $j\in I^+$ there exists $D_{h,j}\subseteq A_j$,
satisfying:
\begin{itemize}
\item $|D_{h,j}|\le 1/(8k^2\tau)$; and
\item $D_{h,j}$ covers $Q_h$;
\item every vertex in $D_{h,j}$ has at most $4k^2\tau|Q_i|$ neighbours in $Q_i$ for all $i\in H\setminus \{h\}$.
\end{itemize}
(This is possible since setting $H=\emptyset$ satisfies the bullets.)
Suppose that there exists
$g\in I^-\setminus H$.

Let $j\in I^+$. For each $h\in H$, each vertex in $Q_h$ has at most $2k\tau|A_j|$ 
neighbours in $A_j$, and so there are at most $2k\tau|Q_h|\cdot|A_j|$ edges between $Q_h$ and $A_j$; and hence at most $|A_j|/(2k)$
vertices in $A_j$ have at least $4k^2\tau|Q_h|$ neighbours in $Q_h$. Let $S_j$ be the set
of vertices $v\in A_j$ such that for each $h\in H$, $v$ has fewer than $4k^2\tau|Q_h|$ neighbours in $Q_h$. It follows that
$|S_j|\ge |A_j|/2$. 

For each $h\in H$ and $j\in I^+$, the set $D_{h,j}$ has cardinality at most $1/(8k^2\tau)$, and since each of its vertices has at most
$\tau|A_g|$ neighbours in $A_g$, it follows that at most $|A_g|/(8k^2)$ vertices in $A_g$ have a neighbour in $D_{h,j}$.
Consequently at most $|A_g|/8$ vertices in $A_g$ have a neighbour in some $D_{h,j}$; and since $|P_g|\ge |A_g|/2$,
there is a subset $T_g\subseteq P_g$ with $|T_g|\ge |A_g|/4$ such that, 
for all $h\in H$ and $j\in J^+$, $T_g$ is anticomplete to $D_{h,j}$.
\\
\\
(1) {\em Let $j\in I^+$. Then there exists $Y_j\subseteq T_g$ with $|T_g\setminus Y_j|< \mu|A_g|$, and $X\subseteq S_j$, such that $X$
covers $Y_j$, and $|X|\le  2 |V_1(G)|^{\phi}/\tau$.}
\\
\\
Choose
$X\subseteq S_j$ maximal such that
\begin{itemize}
\item $|X|\le 2 |V_1(G)|^{\phi}/\tau$; and
\item $|Y|\ge \tau|V_1(G)|^{-\phi}|X|\cdot |A_g|$, where $Y$ is the set of vertices in $T_g$ that have a neighbour in $X$.
\end{itemize}
Suppose that $|T_g\setminus Y|\ge \mu|A_g|$. Since $|S_j|\ge \mu|A_j|$ and
$\tau$ is a $(\phi,\mu)$-band for $(A_i:i\in I)$, it follows that
the max-degree from $S_j$ to $T_g\setminus Y$ is more than $\tau|V_1(G)|^{-\phi}|A_g|$. Choose $v\in S_j$ with more than
$\tau|V_1(G)|^{-\phi}|A_g|$ neighbours in $T_g\setminus Y$. Since $v$ has a neighbour in $T_g\setminus Y$, it follows that $v\notin X$,
and from the maximality of $X$, adding $v$ to $X$ contradicts one of the two bullets in the definition of $X$. The second bullet is
satisfied, and so the first is violated; and hence $|X|+1> 2 |V_1(G)|^{\phi}/\tau$. Since $2 |V_1(G)|^{\phi}/\tau\ge 1$, it follows that $X\ne \emptyset$, and so
$2|X|\ge |X|+1>2 |V_1(G)|^{\phi}/\tau$, and therefore $|X|> |V_1(G)|^{\phi}/\tau$. So $|Y|> \tau|V_1(G)|^{-\phi}(|V_1(G)|^{\phi}/\tau)|A_g|= |A_g|$, a contradiction.
This proves that $|T_g\setminus Y|< \mu|A_g|$, and so proves (1).
\\
\\
(2) {\em There exists $Q_g\subseteq T_g$ with $|Q_g|\ge (64k^2 |V_1(G)|^{\phi})^{-k}|A_g|/8$, and for each $j\in I^+$ there exists a subset
$D_{g,j}\subseteq S_j$ with $|D_{g,j}|\le 1/(8k^2\tau)$, such that $D_{g,j}$ covers $Q_g$.}
\\
\\
For each $j\in I^+$, let $Y_j$ be as in (1), and choose $X_j\subseteq D_j$, such that $X_j$
covers $Y_j$, and $|X_j|\le  2 |V_1(G)|^{\phi}/\tau$. Let $Y$ be the intersection of the sets $Y_j\;(j\in I^+)$. Since each
$Y_j$ satisfies $|T_g\setminus Y_j|< \mu|A_g|$, it follows
that $|Y|\ge |T_g|- k\mu|A_g|\ge |A_g|/8$ (since $k\mu\le 1/8$ and $|T_g|\ge |A_g|/4$).
Let $j\in I^+$. Since $\lfloor 1/(8k^2\tau)\rfloor\ge 1/(16k^2\tau)$ (because $8k^2\tau\le 1$), there is a
partition of $X_j$ into at most $\lceil 16k^2\tau|X_j|\rceil $ sets each of cardinality at most $ 1/(8k^2\tau)$.
But $|X_j|\le  2 |V_1(G)|^{\phi}/\tau$, and so
$$\lceil 16k^2\tau|X_j|\rceil \le \lceil 32k^2 |V_1(G)|^{\phi}\rceil\le 64k^2 |V_1(G)|^{\phi}$$
since $ 32k^2 |V_1(G)|^{\phi}\ge 1$. Thus $X_j$ admits a partition $\mathcal{R}_j$ into at most $64k^2 |V_1(G)|^{\phi}$
sets each of cardinality at most  $ 1/(8k^2\tau)$.
For each $v\in Y$, there exists $u\in X_j$ adjacent to $Y$; choose some such $u$, choose $R\in \mathcal{R}_j$
containing $u$, and say $R$ is the {\em $j$-type} of $v$. Each vertex of $v$ has a $j$-type, for each $j\in I^+$; and since
there are only at most $64k^2 |V_1(G)|^{\phi}$ $j$-types for each $j$, and $|I^+|\le k$,
it follows that there exists $Q_g\subseteq Y$ with
$$|Q_g|\ge (64k^2 |V_1(G)|^{\phi})^{-k}|Y| \ge (64k^2 |V_1(G)|^{\phi})^{-k}|A_g|/8,$$
such that for all $j\in I^+$, all members of $Q_g$ have the same $j$-type, say $D_{g,j}\in \mathcal{R}_j$,
and each $D_{g,j}$ covers $B_g$.
This proves (2).

\bigskip
From (2), this contradicts the maximality of $H$. (Note that for $j\in I^+$ and $h\in H$, every vertex in $D_{h,j}$
has no neighbours in $Q_g$, since $Q_g\subseteq T_g$; and every vertex of $D_{g,j}$ has at most $4k^2\tau|Q_h|$ neighbours in $Q_h$,
because $D_{g,j}\subseteq S_j$.) This proves that $H=I^-$. Thus we have shown that
for all $h\in I^-$ there exists $Q_h\subseteq P_h$ with
$|Q_h|\ge (64k^2 |V_1(G)|^{\phi})^{-k}|A_h|/8$, and for all $h\in I^-$ and $j\in I^+$ there exists $D_{h,j}\subseteq A_j$,
satisfying:
\begin{itemize}
\item $|D_{h,j}|\le 1/(8k^2\tau)$;
\item $D_{h,j}$ covers $Q_h$; and
\item every vertex in $D_{h,j}$ has at most $4k^2\tau|Q_i|$ neighbours in $Q_i$ for all $i\in H\setminus \{h\}$.
\end{itemize}
Now let $i\in I^-$.
Since there are only at most $k^2$ sets $D_{h,j}$, and each has cardinality at most $1/(8k^2\tau)$, and every
vertex in a set $D_{h,j}$ for $h\ne i$ has at most $4k^2\tau|Q_i|$ neighbours in $Q_i$, it follows that at most $|Q_i|/2$
vertices in $Q_i$ have a neighbour in some set $D_{h,j}$ with $h\ne i$; and so there exists $C_i\subseteq Q_i$ with
$|C_i|\ge |Q_i|/2$ such that $C_i$ is anticomplete to $D_{h,j}$ for all $h\in I^-\setminus \{i\}$ and $j\in I^+$.
Hence $|C_i|\ge (64k^2)^{-k}|V_1(G)|^{-k\phi}|A_h|/16$ for each $i\in I^-$.
Moreover, since the union of all the sets $D_{h,j}$ has cardinality at most $1/(8\tau)$, and each vertex of this union has at most $\tau|A_i|$
neighbours in $A_i$, it follows that at most $|A_i|/8$ vertices in $A_i$ have a neighbour that belongs to some $D_{h,j}$;
and consequently there exists $B_i\subseteq A_i$ with $C_i\subseteq B_i$, and with $|B_i|\ge |A_i|/2$, such that $B_i\setminus C_i$
is anticomplete to all the sets $D_{h,j}$.
This proves \ref{newleafcover}.~\bbox

By combining \ref{homog} and \ref{newleafcover}, and fixing values for $\phi$ and $\mu$, we obtain:
\begin{thm}\label{combined}
Let $k\ge 0$ be an integer, and let  $c>0$.
Then there exists an integer $K>0$ with the following property.
Let $\mathcal{A}=(A_i:i\in I)$ be a
parade  of length at least $(K,K)$ in a bigraph $G$. Let $\beta=(8k)^{-1-2K^2k/c}$.
Then either
\begin{itemize}
\item there exist $h\in I^-$ and $j\in I^+$, and $X\subseteq A_h$ and $Y\subseteq A_j$ with
$\frac{|X|}{|A_h|}, \frac{|Y|}{|A_j|}\ge \beta$, such that $X,Y$ are anticomplete; or
\item there exist $h\in I^-$ and $j\in I^+$ such that some $v\in A_j$ has at least $\frac{\beta}{8k^2}|A_h|$ neighbours in $A_h$; or
\item there exist $J\subseteq I$ with $|J^-|=|J^+|=k$, 
and for each $h\in J^-$ there exists $B_h\subseteq A_h$ with $|B_h|\ge \beta |A_h|/2$, and there exists
$C_h\subseteq B_h$ with $|C_h|\ge \beta|V_1(G)|^{-c}|A_h|/16$; and
for each $h\in J^-$ and $j\in J^+$ there exists $D_{h,j}\subseteq A_j$ covering $C_h$, such that
$D_{h,j}$ is anticomplete to $B_i\setminus C_i$ for all $i\in J^-$, and is anticomplete to $C_i$
for all $i\in J^-\setminus \{h\}$.
\end{itemize}
\end{thm}
\Proof
By \ref{homog}, taking $\mu=1/(8k)$ and $\phi=c/k$, we may assume that there exist $J\subseteq I$ with $|J^-|=|J^+|=k$, and a subset $F_i\subseteq A_i$ with $|F_i|\ge \beta |A_i|$ for each $i\in J$, 
such that
$(F_i:i\in J)$ has a $(\phi,\mu)$-band $\tau$. We may assume that 
for all $j\in J^+$ and $v\in A_j$ and $h\in I^-$, $v$ has fewer than $(\beta/(8k^2))|A_h|$ neighbours in $A_h$, and hence has
fewer than $|F_h|/(8k^2)$ neighbours in $F_h$.
Consequently we may assume that $\tau\le 1/(8k^2)$.
By \ref{newleafcover} applied to $\mathcal{F}=(F_i:i\in J)$, 
for each $h\in J^-$ there exists $B_h\subseteq F_h$ with $|B_h|\ge |F_h|/2\ge \beta|A_h|/2$, and there exists
$C_h\subseteq B_h$ with 
$$|C_h|\ge |V_1(G)|^{-k\phi}|F_h|/16\ge \beta |V_1(G)|^{-k\phi}|A_h|/16;$$ 
and
for each $h\in J^-$ and $j\in J^+$ there exists $D_{h,j}\subseteq F_j$ covering $C_h$, such that
$D_{h,j}$ is anticomplete to $B_i\setminus C_i$ for all $i\in I^-$, and is anticomplete to $C_i$
for all $i\in J^-\setminus \{h\}$. Then
the theorem is satisfied.
This proves \ref{combined}.~\bbox

\section{The proof of \ref{linmainthm}}

If $T$ is a tree and $w\in V(T)$,
we say that the {\em $w$-radius} of $T$ is the maximum integer $r$ such that some path of $T$ with one end $w$ has $r$ edges.
For $v\in V(T)\setminus \{w\}$, the {\em $w$-parent} of $v$ is the neighbour of $v$ in              
the path of $T$ between $v,w$. We define $d_T(u,v)$ to be the distance in $T$ between $u,v$.

If $G$ is a bigraph and $T$ is an induced sub-bigraph that is a tree bigraph, we say that $T$ is a {\em induced subtree} of $G$.
If $\mathcal{B}$ is a parade in a bigraph $G$, an induced subtree $T$ of $G$ is {\em $\mathcal{B}$-rainbow} if
every vertex of $H$ belongs to some block of $\mathcal{B}$, and every block of $\mathcal{B}$ contains at most one vertex of $H$.

If $I$ is a set of nonzero integers, a {\em shape} in $I$ is a tree $S$ with $V(S)\subseteq I$, such that for every
edge $ij$ of $S$, $i$ and $j$ have opposite sign.
Let $\mathcal{A}=(A_i:i\in I)$ be a parade in a bigraph $G$,
let $T$ be an $\mathcal{A}$-rainbow induced subtree of $G$, and let $S$ be a shape in $I$. We say that $S$ is the {\em shape of $T$} if
\begin{itemize}
\item for each $i\in I$, $i\in V(S)$ if and only if some vertex of $A_i$ belongs to $V(T)$; and
\item for all $i,j\in I$, $i$ is adjacent to $j$ in $S$ if and only if there exists $u\in A_i$ and $v\in A_j$ such that
$uv$ is an edge of $T$.
\end{itemize}
Thus every induced subtree of $G$ has a unique shape.

Let $\mathcal{A}^0=(A^0_i:i\in I)$ be a parade in a bigraph $G$, and let $r\ge 0$ be an integer.
For each $i \in I$, and for $1\le q\le r$ let $A^{q}_i\subseteq A^{q-1}_i$ be nonempty. For $0\le q\le r$
let $\mathcal{A}^q$ be the parade $(A^q_i:i\in I)$. We call the sequence $(\mathcal{A}^0\ll \mathcal{A}^r)$ a {\em nested parade sequence}.
Let $(\mathcal{A}^0\ll \mathcal{A}^r)$ be a nested parade sequence with notation as above, let $T$ be an $\mathcal{A}^0$-rainbow 
induced subtree of $G$, and let
$w\in V(T)$, such that $T$ has $w$-radius at most $r$. We say that $T$ is {\em $w$-isolated in $(\mathcal{A}^0\ll \mathcal{A}^r)$} if:
\begin{itemize}
\item for each $v\in V(T)$, $v\in A_i^{r-d_T(v,w)}$ for some $i\in I$; and 
\item for each $v\in V(T)\setminus \{w\}$, let $q=r+1-d_T(v,w)$; for $i\in I$, $v$ has a neighbour in $A^{q}_i$ only if
the $w$-parent of $v$ in $T$ also belongs to $A^{q}_i$.
\end{itemize}

Again, let  $(\mathcal{A}^0\ll \mathcal{A}^r)$ be a nested parade sequence with notation as before. Let $i\in I$; we say that
a vertex $w\in A^0_i$ is {\em $r$-panarboreal in  $(\mathcal{A}^0\ll \mathcal{A}^r)$} 
if for every shape $S$ in $I$ with $i\in V(S)$
and with $i$-radius at most $r$,
there is an $\mathcal{A}^0$-rainbow induced subtree $T$ of $G$ with shape $S$ that is $w$-isolated in  
$(\mathcal{A}^0\ll \mathcal{A}^r)$.
We will prove:

\begin{thm}\label{rainbow}
Let $0<c\le 1$. For all integers $r,k\ge 0$, 
there exist an integer $K>0$, and $\gamma>0$ with the following property. 
Let $G$ be a bigraph and let $\mathcal{A}=(A_i:i\in I)$ be a parade in $G$ with length at least $(K,K)$. 
Then either:
\begin{itemize}
\item there exist $h\in I^-$ and $j\in I^+$, and $X\subseteq A_h$ and $Y\subseteq A_j$, either with
$|X|\ge \gamma |A_h|$ and $|Y|\ge \gamma|V_2(G)|^{-c}|A_j|$, or with 
$|X|\ge \gamma |V_1(G)|^{-c}|A_h|$ and $|Y|\ge \gamma|A_j|$, such that $X,Y$ are anticomplete; or
\item there exist $h,j\in I$ with opposite sign, and $v\in A_h$, such that $v$ has at least $\gamma|A_j|$ neighbours in~$A_j$; or
\item there exist $J\subseteq I$ with $|J^-|,|J^+|\ge k$, and a nested parade sequence  $(\mathcal{A}^0\ll \mathcal{A}^r)$ with
$\mathcal{A}^q=(A^q_i:j\in J)$ for $0\le q\le r$, and with the following properties.
For each $j\in J$, $A^0_j\subseteq A_j$ and $|A^r_j|\ge \gamma|A_j|$;
and for each $h\in J^-$ there exists $C_h\subseteq A^r_h$ with 
$|C_h|\ge \gamma |V_1(G)|^{-c} |A_h|$, such that every vertex in $C_h$ is $r$-panarboreal in $(\mathcal{A}^0\ll \mathcal{A}^r)$.
\end{itemize}
\end{thm}










\Proof
We will prove, by induction on $r$, that for each value of $r$ the statement holds for all $k$. 
If $r=0$ then the third outcome of the theorem is true, setting $K=k$, and $\gamma=1$,
and taking $J=I$, and $A^0_h=A_h$ for all $h\in I$, and $C_h=A_h$ for all $h\in J^-$. Thus we may assume that $r\ge 1$, 
and the claim holds for $r-1$ and all $k$. 

Choose $k'$ such that setting $K=k'$ satisfies \ref{combined}. From the inductive hypothesis,
there exist an integer $K>0$ 
and $\gamma'>0$ such that the assertion of \ref{rainbow} holds with $r,k, \gamma$ 
replaced by $r-1,k', \gamma'$ respectively. 
Let $\beta=(8k)^{-1-2k'^2k/c}$, and 
$\gamma= \frac{\beta}{8k^2}\gamma'$.
We claim that the theorem is satisfied.
To see this, let $G$ be a bigraph and let $\mathcal{A}=(A_i:i\in I)$ be a parade in $G$ with length at least $(K,K)$. 
\\
\\
(1) {\em We may assume that there exist $L\subseteq I$ with $|L^-|,|L^+|\ge k'$, and a nested parade sequence 
$(\mathcal{B}^0\ll \mathcal{B}^{r-1})$, with  $\mathcal{B}^q=(B^q_i:j\in L)$ for $0\le q\le r-1$, and with the following properties.
For each $j\in L$, $B^0_j\subseteq A_j$ and $|B^{r-1}_j|\ge \gamma|A_j|$;
and for each $h\in L^-$ there exists $C_h\subseteq B^{r-1}_h$ with          
$|C_h|\ge \gamma |V_1(G)|^{-c} |A_h|$, such that every vertex in $C_h$ is $(r-1)$-panarboreal in $(\mathcal{B}^0\ll \mathcal{B}^{r-1})$.}
\\
\\
Let
$G^T$ be the bigraph obtained from $G$ by setting $V_1(G^T)=V_2(G)$ and $V_2(G^T)=V_1(G)$, and let $I^T$ be the set of all integers
$-i$ where $i\in I$. Thus $\mathcal{A}^T=(A_i:i\in I^T)$ is a parade in $G^T$. 
From the choice of $K$, applied to $G^T$ and $\mathcal{A}^T$, either 
\begin{itemize}
\item there exist $h\in I^-$ and $j\in I^+$, and $X\subseteq A_h$ and $Y\subseteq A_j$, either with
$|X|\ge \gamma' |A_h|$ and $|Y|\ge \gamma'|V_2(G)|^{-c}|A_j|$, or with
$|X|\ge \gamma' |V_1(G)|^{-c}|A_h|$ and $|Y|\ge \gamma'|A_j|$, 
such that $X,Y$ are anticomplete; or
\item there exist $h,j\in I$ with opposite sign, and $v\in A_h$, such that $v$ has at least $\gamma'|A_j|$ neighbours in $A_j$; or
\item the statement of (1) holds.
\end{itemize}
In the first case, since $\gamma\le \gamma'$, the first outcome of the theorem holds, and similarly
in the second case, the second outcome of the theorem holds. So we may assume that the third case holds.
This proves (1).
\\
\\
(2) {\em We may assume that there exists $J\subseteq L$ with $|J^-|=|J^+|=k$,
and for each $h\in J^-$ there exists $B^r_h\subseteq B^{r-1}_h$ with $|B^r_h|\ge \beta |B^{r-1}_h|/2$, and there exists
$C_h\subseteq B^r_h$ with $|C_h|\ge \beta|V_1(G)|^{-c}|B^{r-1}_h|/16$; and
for each $h\in J^-$ and $j\in J^+$ there exists $D_{h,j}\subseteq C_j$ covering $C_h$, such that
$D_{h,j}$ is anticomplete to $B^r_i$ for all $i\in J^-\setminus \{h\}$.}
\\
\\
Let $A_i'=B^{r-1}_i$ for $i\in L^-$, and $A_i'=C_i$ for $i\in L^+$. Then $(A_i':i\in L)$ is a parade of length at least $(k',k')$,
and so from the choice of $k'$, either
\begin{itemize}
\item there exist $h\in L^-$ and $j\in L^+$, and $X\subseteq A_h'$ and $Y\subseteq A_j'$ with
$\frac{|X|}{|A_h'|}, \frac{|Y|}{|A_j'|}\ge \beta$, such that $X,Y$ are anticomplete; or
\item there exist $h\in L^-$ and $j\in L^+$ such that some $v\in A_j$ has at least $\frac{\beta}{8k^2}|A_h'|$ neighbours in $A_h'$; or
\item the statement of (2) holds.
\end{itemize}
In the first case, the first outcome of the theorem holds, since 
$$|X|\ge \beta |A_h'|=\beta|B^{r-1}_h| \ge \beta \gamma'|A_h| \ge \gamma |A_h|$$
(because $\beta \gamma'\ge \gamma$), and 
$$|Y|\ge \beta |A_j'|=\beta |C_j|\ge \beta \gamma' |V_2(G)|^{-c} |A_j|\ge \gamma|V_2(G)|^{-c}|A_j|.$$
In the second case, the second outcome of the theorem holds, since
$$\frac{\beta}{8k^2}|A_h'|= \frac{\beta}{8k^2}|B^{r-1}_h|\ge \frac{\beta}{8k^2}\gamma'|A_h|\ge \gamma |A_h|.$$
Thus we may assume that the third case holds. This proves (2).

\bigskip
Define $B^r_j=B^{r-1}_j$ for each $j\in J^+$. Then $|B^r_i|\ge \gamma|A_i|$ for each $i\in J$, since for $h\in J^-$, 
$$|B^r_h|\ge  \beta |B^{r-1}_h|/2\ge (\beta/2)\gamma'|A_h|\ge \gamma|A_h|$$
and for each $j\in J^+$, 
$$|B^r_j|=|B^{r-1}_j|\ge \gamma'|A_j|\ge \gamma|A_j|.$$
Moreover, for each $h\in J^-$, $C_h\subseteq B^r_h$, and 
$$|C_h|\ge \beta|V_1(G)|^{-c}|B^{r-1}_h|/16\ge  \beta|V_1(G)|^{-c}\gamma'|A_h|/16\ge \gamma |V_1(G)|^{-c} |A_h|.$$
Let $\mathcal{A}^q=(B^q_j:j\in J)$ for $0\le q\le r$.
To complete the proof, we will show that for each $h\in J^-$ and each $w\in C_h$, $w$ is $r$-panarboreal in $(\mathcal{A}^0\ll \mathcal{A}^{r})$.

Let $S$ be a shape in $J$ with $h\in V(S)$ and with $h$-radius at most $r$. Let $j_1\ll j_t$ be the neighbours of $h$
in $S$; thus $j_1\ll j_t\in J^+$. For $1\le s\le t$, let $S_s$ be the component of $S\setminus \{h\}$ that contains $j_s$.
Thus $S_s$ has $j_s$-radius at most $r-1$, and $S_s$ is a shape in $J$ and hence in $L$. 
For $1\le s\le t$, since $D_{h,j_s}$
covers $C_h$, there is a vertex $w_s\in D_{h,j_s}$ adjacent to $w$. Consequently $w_s$ has no neighbours in 
$B^r_i\setminus C_i$ for all $i\in J^-$, and has no neighbours in $C_i$
for all $i\in J^-\setminus \{h\}$. Since every vertex in $C_i$ is $(r-1)$-panarboreal in $(\mathcal{B}^0\ll \mathcal{B}^{r-1})$,
it follows that 
there is an $\mathcal{A}$-rainbow subtree $T_s$ of $G$ with shape $S_s$ and with $w_s\in V(T_s)$, such that
$T_s$ is $w_s$-isolated in $(\mathcal{B}^0\ll \mathcal{B}^{r-1})$. 
Let $T$ be the tree obtained from the union of the trees $T_1\ll T_t$
by adding the vertex $w$ and the edges $ww_1\ll ww_s$. 
\\
\\
(3) {\em $T$ is an induced subtree of $G$.}
\\
\\
To see this, since 
each $T_s$ is induced,  it suffices to check
that 
\begin{itemize}
\item $w$ has a unique neighbour $w_s$ in $V(T_s)$, for $1\le s\le t$; and 
\item there is no edge of $G$ between 
$V(T_s)$ and $V(T_{s'})$ for distinct $s,s'$. 
\end{itemize}
To prove the first, suppose that $w$ is adjacent to some vertex $v\in V(T_s)$ where $v\ne w_s$.
Since $T_s$ is $w_s$-isolated in $(\mathcal{B}^0\ll \mathcal{B}^{r-1})$, and $w\in B^{r-1}_h$, it follows that
the $w_s$-parent of $v$ in $T_s$ also belongs to $B^{r-1}_h$; but this is not the case since $h\notin V(S_s)$. So 
 $w$ has a unique neighbour $w_s$ in $V(T_s)$, for $1\le s\le t$. This proves the first bullet.
Now suppose that $1\le s,s'\le t$ with $s\ne s'$, and some vertex
$u\in V(T_s)$ is adjacent to some vertex $v$ in $V(T_{s'})$. Since $G$ is bipartite, one of $u,v$ is closer in $T$ to $w$ than 
the other, say $d_T(v,w)<d_T(u,w)$, and consequently $u\ne w_s$.
Since $T_s$ is $w_s$-isolated in $(\mathcal{B}^0\ll \mathcal{B}^{r-1})$, and $u\ne w_s$, it follows that $v$ and the 
$w_s$-parent of $u$ in $T_s$
belong to the same block of $\mathcal{A}$, a contradiction since the shapes of $T_s,T_{s'}$ are vertex-disjoint. This proves the
second bullet above, and so proves (3).
\\
\\
(4) {\em $T$ is $w$-isolated in $(\mathcal{A}^0\ll \mathcal{A}^{r})$.}
\\
\\
We must show that:
\begin{itemize}
\item for each $v\in V(T)$, $v\in B_i^{r-d_T(v,w)}$ for some $i\in J$; and
\item for each $v\in V(T)\setminus \{w\}$, let $q=r+1-d_T(v,w)$; for $i\in J$, $v$ has a neighbour in $B^{q}_i$ only if
the $w$-parent of $v$ in $T$ also belongs to $B^{q}_i$.
\end{itemize}

For the first bullet, 
since $w\in C_h\subseteq B^r_h$, we may assume that $v\ne w$; let $v\in V(T_s)$ where $1\le s\le t$. 
Then the first bullet follows, since
$T_s$ is $w_s$-isolated in $(\mathcal{B}^0\ll \mathcal{B}^{r-1})$, and 
$$r-d_T(v,w) =(r-1)-d_{T_s}(v,w_s).$$

For the second bullet, let $v\in V(T_s)$.
Since, as we saw earlier, $w_s$ has no neighbours in
$B^r_j$ for all $j\in J^-\setminus \{h\}$, we may assume that $v\ne w_s$.
But then the second bullet is true,
since
$T_s$ is $w_s$-isolated in $(\mathcal{B}^0\ll \mathcal{B}^{r-1})$, and 
$$r-d_T(v,w) =(r-1)-d_{T_s}(v,w_s),$$ 
and 
the $w_s$-parent of $v$ in $T_s$ is the $w$-parent of $v$ in $T$.
This proves (4).

\bigskip

From (3) and (4), this 
proves \ref{rainbow}.~\bbox

Now we can prove \ref{linmainthm}, which we restate:
\begin{thm}\label{linmainthm2}
Let $T$ be an ordered tree bigraph. For all  $c>0$ there exists $\vare>0$ with the following property.
Let $G$ be an ordered bigraph not containing $T$, such that every vertex in $V_1(G)$ has degree less than $\vare |V_2(G)|$, and 
every vertex in $V_2(G)$ has degree less than $\vare |V_1(G)|$.
Then
there are subsets $Z_i\subseteq V_i(G)$ for $i = 1,2$, either with $|Z_1|\ge \vare |V_1(G)|$ and $|Z_2|\ge \vare |V_2(G)|^{1-c}$, 
or with $|Z_1|\ge \vare |V_1(G)|^{1-c}$ and $|Z_2|\ge \vare |V_2(G)|$, such that $Z_1,Z_2$ are anticomplete.
\end{thm}
\Proof We may assume that $|V(T)|\ge 2$; choose an integer $r$ such that $T$ has $w_1$-radius at most $r$, for some vertex $w_1\in V_1(T)$,
and  $T$ has $w_2$-radius at most $r$, for some vertex $w_2\in V_2(T)$. 
Choose an integer $k$ 
such that $|V_1(T)|, |V_2(T)|\le k$. Choose $K, \gamma$ as in \ref{rainbow}. Let $\vare=\gamma/(2K)$. We claim that $\vare$
satisfies the theorem.

Let $G$ be an ordered bigraph that does not contain $T$, such that 
every vertex in $V_1(G)$ has degree less than $\vare |V_2(G)|$, and
every vertex in $V_2(G)$ has degree less than $\vare |V_1(G)|$.
If $G$ has no edges then $Z_1,Z_2$ exist as required, so we may assume that $G$ has an edge; and so $\vare|V_i(G)|<1$ for $i = 1,2$.
Let $p=\lceil |V_1(G)|/(2K)\rceil$; then $p\le |V_1(G)|/K$, since $|V_1(G)|> 1/\vare\ge K$. Let the vertices of $V_1(G)$
be $u_1\ll u_{n_1}$, ordered according to the linear order of $V_1(G)$ imposed by $G$. 
Let 
$$A_{i}=\{u_{(K-i)p+1}\ll u_{(K-i+1)p}\}$$ 
for $-K\le i\le -1$. Similarly, let $q= \lceil |V_2(G)|/(2K)\rceil$, and 
$V_2(G)=\{v_1\ll v_{n_2}\}$ in order, 
and for $1\le i\le K$ let 
$$A_i=\{v_{(i-1)q+1}\ll v_{iq}\}.$$ 
Let $I=\{-K\ll -1,1\ll K\}$. Then $\mathcal{A}=(A_i:i\in I)$ is a parade in $G$,
of length $(K,K)$ and width $(p,q)$, and all its blocks are intervals of the linear order, in the natural sense.
Since the blocks of $(A_i:i\in J)$ are intervals and are numbered in order, it follows that:
\\
\\
(1) {\em Let $R$ be an $\mathcal{A}$-rainbow induced subtree of $G$ with shape $S$. The orders of $V_1(G)$ and $V_2(G)$
induce orders on $V_1(R), V_2(R)$, making $R$ into an ordered tree bigraph $R'$. Also the orders
on $I^-$ and $I^+$ make $S$ into an ordered tree bigraph $S'$; and $R'$ is isomorphic to $S'$.}
\\
\\
Certainly $R$ is isomorphic to $S$, but we need to check that the natural isomorphism preserves the vertex-orders. For
each $v\in V(R)$, let $f(v)\in I$ such that $v\in A_{f(v)}$; then $f$ is an isomorphism from $R$ to $S$. Let 
$u,v\in V_1(R)$ say, where $u$ is earlier than $v$ in the order that $R'$ imposes on $V_1(R')$. Hence 
$u$ is earlier than $v$ in the order that $G$ imposes on $V_1(G)$, that is, $i<j$ where  $u=u_i$ and $v=u_j$. 
Since the blocks of $\mathcal{A}$
are intervals, numbered in order, and $i<j$, and $f(u)\ne f(v)$ since $R$ is $\mathcal{A}$-rainbow, 
it follows that $f(u)<f(v)$. Hence $f(u)$ is earlier than $f(v)$ in the order imposed on $V_1(S')$
by $S'$; and so $f$ maps the order of $V_1(R')$ imposed by $R'$ to the order of $V_1(S')$ imposed by $S'$. Similarly
$f$ maps the order of $V_2(R')$ imposed by $R'$ to the order of $V_2(S')$ imposed by $S'$. This proves (1).

\bigskip

Let us apply \ref{rainbow} to $(A_i:i\in I)$, and deduce that one of the three outcomes of \ref{rainbow} holds.
Suppose that the first outcome holds, that is, there exist $h\in I^-$ and $j\in I^+$, and $Z_1\subseteq A_h$ and $Z_2\subseteq A_j$, 
either with
$|Z_1|\ge \gamma |A_h|$ and $|Z_2|\ge \gamma|V_2(G)|^{-c}|A_j|$, or with 
$|Z_1|\ge \gamma |V_1(G)|^{-c}|A_h|$ and $|Z_2|\ge \gamma|A_j|$, such that $Z_1,Z_2$ are anticomplete, and from the symmetry 
we assume the first.
Since 
$$\gamma |A_h|=\gamma p \ge \gamma |V_1(G)|/(2K) = \vare|V_1(G)|$$
and 
$$\gamma|V_2(G)|^{-c}|A_j|\ge \gamma|V_2(G)|^{-c}q\ge \gamma|V_2(G)|^{-c} |V_2(G)|/(2K)= \vare|V_2(G)|^{1-c},$$ 
in this case the theorem holds.

Now suppose that the second outcome holds, that is, there exist $h,j\in I$ with opposite sign, and $v\in A_h$, such that $v$ 
has at least $\gamma|A_j|$ neighbours in~$A_j$. From the symmetry we may assume that $h\in I^-$. Since 
$$\gamma|A_j|=\gamma q\ge \gamma |V_2(G)|/(2K)= \vare |V_2(G)|$$
this is impossible.

Finally, suppose that the third outcome holds, that is, there exist $J\subseteq I$ with $|J^-|,|J^+|\ge k$, 
and a nested parade sequence  $(\mathcal{A}^0\ll \mathcal{A}^r)$ with $\mathcal{A}^q=(A^q_i:j\in J)$ for $0\le q\le r$, 
and with the following properties.
For each $j\in J$, $A^0_j\subseteq A_j$ and $|A^r_j|\ge \gamma|A_j|$;
and for each $h\in J^-$ there exists $C_h\subseteq A^r_h$ with
$|C_h|\ge \gamma |V_1(G)|^{-c} |A_h|$, such that every vertex in $C_h$ is $r$-panarboreal in $(\mathcal{A}^0\ll \mathcal{A}^r)$.
Since $|V_1(T)|, |V_2(T)|\le k$, it follows that there is a shape $S$ in $J$, such that the ordered tree bigraph $T$ is isomorphic 
to the ordered tree bigraph
$S'$ obtained from $S$ as in (1). Let this isomorphism map $w_1$ to $h\in J^-$,
and choose $w\in C_h$. Since $w$ is $r$-panarboreal in $(\mathcal{A}^0\ll \mathcal{A}^r)$, it follows that there is an 
$(A_i:i\in I)$-rainbow induced subtree $R$ of $G$, with shape $S$. Let $R'$ be as in (1); then $R'$ is isomorphic to $S'$
by (1), and hence isomorphic to $T$, and therefore $G$ contains $T$, a contradiction.
This proves~\ref{linmainthm2}.~\bbox

\end{document}